\documentclass[11pt,amsfonts]{article}
\usepackage{graphicx}
\usepackage{latexsym}
\usepackage{amssymb} 
\usepackage{amsmath}
\usepackage{layout}  
\newtheorem{prop}{Proposition}
\newtheorem{lemma}{Lemma} 
\newtheorem{definition}{Definition}

\newtheorem{theorem}{Theorem}

\newtheorem{remark}{Remark}

 \newcommand{\fin}
{ \vspace{-0.6cm}
\begin{flushright}
\mbox{$\Box$}
\end{flushright}
\noindent }

\newcommand{\R}{\mathbb{R}}
\newcommand{\N}{\mathbb{N}}

\newcommand{\DE}{\widehat{A}^{n}}
\newcommand{\XE}{\overline{X}^{n}}

\newcommand{\eps}{\varepsilon}
\def\D{\Delta}
\def\e{{\cal \varepsilon}}
\def\EE{{\cal E}}

\def\HH{{\cal H}}

\def\N{{\mathbb{N}}}

\def\R{{\mathbb{R}}}

\def\lcr{\left[}
\def\lpa{\langle}
\def\lva{\left|}

\def\rcr{\right]}
\def\rpa{\rangle}
\def\rva{\right|}

\def\sk{{\mathbb{D}}}
\def\Un{{\bf 1}}

\newcommand{\ignore}[1]{}
\begin{document}

\renewcommand{\thefootnote}{\fnsymbol{footnote}}

\title{{\bf Exact rate of convergence of some approximation schemes associated
to SDEs driven by a fractional Brownian motion}}
\date{
\small {\bf Andreas Neuenkirch}\\
Johann Wolfgang Goethe Universit\"at, Fachbereich Mathematik, \\
Robert-Mayer-Stra{\ss}e 10, 60325 Frankfurt am Main, Germany\\
{\tt neuenkirch@math.uni-frankfurt.de}\\
$\mbox{ }$\\ 
$\mbox{ }$\\
\small {\bf Ivan Nourdin}\\
LPMA, Universit\'e Pierre et Marie Curie Paris 6,\\
Bo\^ite courrier 188, 4 Place Jussieu, 75252 Paris Cedex 5, France\\
{\tt nourdin@ccr.jussieu.fr}\\
}

 \maketitle

\begin{abstract}
\noindent
In this paper, we derive the exact rate of convergence of some approximation schemes associated to scalar stochastic differential 
equations driven by a fractional Brownian motion with Hurst index $H$. We consider two cases. If $H>1/2$, 
the exact rate of convergence of the Euler scheme is determined. We show that  the error of the Euler scheme  converges almost surely 
to a random variable, which in particular depends on the Malliavin derivative of the solution.
This result extends  those contained in \cite{neuenkirch} and \cite{cras}.
When $1/6<H<1/2$, the exact rate of convergence of the Crank-Nicholson scheme is determined for a particular equation. Here we show convergence in law of the error to a random variable, which depends on the solution of the equation and an independent Gaussian random variable.
\\
\end{abstract}



\eject

\section{Introduction}

Let $B=(B_t,\,t\in [0,1])$ be a fractional Brownian motion (in short: fBm) with Hurst parameter $H\in (0,1)$, i.e., $B$ is  a continuous centered Gaussian process with covariance function
$$
R_{H}(s,t)=\frac{1}{2}(s^{2H}+t^{2H}-|t-s|^{2H}),\quad s,t\in [0,1].
$$
For $H=1/2$, $B$ is a standard Brownian motion, while for $H\neq 1/2$, it is neither a semimartingale
nor a Markov process. Moreover, it holds
$$
( { \rm E} |B_t-B_s|^2)^{1/2}=|t-s|^H,\quad s,t\in [0,1],
$$
and almost all sample paths of $B$ are H\"older continuous of any order $\alpha\in(0,H)$.

In this paper, we are interested in the pathwise approximation of the equation
\begin{equation}\label{eq1}
X_t=x_0+\int_0^t \sigma(X_s)dB_s +\int_0^t b(X_s)ds,\quad t\in [0,1],
\end{equation}
with a deterministic initial value $x_0\in\R$. Here, $\sigma$ and $b$ satisfy some standard smoothness
assumptions and the integral equation (\ref{eq1}) is understood in the sense of 
Russo-Vallois. Let us recall briefly the significant points of this theory.

\medskip

\begin{definition}(following \cite{RV93})
Let $Z=(Z_t)_{t\in [0,1]}$ be a stochastic process with continuous paths. 
\begin{description}
\item $\bullet$ A family of processes $(H_t^{(\e)})_{t\in[0,1]}$ is said to converge to the process $(H_t)_{t\in [0,1]}$ in the
{\rm ucp sense}, if $\sup_{t \in [0,1]} |H_t^{(\e)}-H_t|$ goes to 0 in probability, as $\e\rightarrow 0$.
\item $\bullet$
The {\rm (Russo-Vallois) forward integral} $\int_0^t Z_sd^- B_s$ is defined by
\begin{equation}\label{forward}
\lim_{\e\rightarrow 0}-{\rm ucp}\,\,\,\e^{-1}\int_0^t Z_t(B_{t+\e}-B_t)dt,
\end{equation}
provided the limit exists.

\item $\bullet$
The {\rm (Russo-Vallois) symmetric integral} $\int_0^t Z_sd^\circ B_s$ is defined by
\begin{equation}\label{symetric}
\lim_{\e\rightarrow 0}-{\rm ucp}\,\,\,(2\e)^{-1}\int_0^t (Z_{t+\e}+Z_t)(B_{t+\e}-B_t)dt,
\end{equation}
provided the limit exists.
\end{description}
\end{definition} 
\medskip
Now we state the exact meaning of  equation (\ref{eq1}) and  give conditions for the existence and uniqueness of its solution.
We consider two cases, according to the value of $H$:
\begin{itemize}
\item {\it Case} $H>1/2$.
\smallskip

Here the integral with respect
to $B$ is defined by the forward integral (\ref{forward}).
\begin{prop}\label{prop17}
If $\sigma\in \mathcal{C}^2_b$ and if $b$ satisfies a global Lipschitz condition, then  the equation 
\begin{equation}\label{eq1f}
X_t=x_0+\int_0^t \sigma(X_s)d^- B_s + \int_{0}^{t}b(X_{s}) \, ds, \quad t\in [0,1]
\end{equation}
 admits a unique
solution $X$ in the set of processes whose paths are H\"older continuous of order $\alpha>1-H$. 
Moreover, we have a Doss-Sussmann type \cite{D,S}
representation:
\begin{equation}\label{ds}
X_t=\phi(A_t,B_t),\quad t\in [0,1],
\end{equation}
where  $\phi$ and $A$ are given respectively by
\begin{equation}\label{phi}
\frac{ \partial  \phi}{ \partial x_{2}}(x_{1},x_{2})=\sigma(\phi(x_{1},x_{2})),\quad \phi(x_{1},0)=x_{1}, \qquad x_{1},x_{2}\in\R
\end{equation}
and
\begin{equation}\label{D}
A'_t={\rm exp}\left(
-\int_0^{B_t} \sigma'(\phi(A_t,s))ds
\right) b(\phi(A_t,B_t)), \quad A_0=x_0, \qquad t \in [0,1].
\end{equation}
\end{prop}
{\bf Proof}. If $X$ and 
$Y$ are two real processes whose paths are  a.s.   H\"older continuous of index $\alpha>0$
and $\beta>0$ with $\alpha+\beta>1$, then $\int_0^t Y_s d^{-}X_s$ coincides with the Young integral
$\int_0^t Y_sdX_s$ (see \cite{RV05}, Proposition 2.12). Consequently, Proposition \ref{prop17} is
a consequence of, e.g., \cite{ZK} or \cite{NR}.\fin

\item {\it Case $1/6<H<1/2$}.
\smallskip

When $H<1/2$, in particular the forward integral $\int_0^t B_sd^-B_s$ does not exist.
Thus, in this case, the use of the symmetric integral (\ref{symetric}) is more adequate. 
Here we consider only the case $b=0$: for the general case see  \cite{lpma1024}, 
\cite{nosi} and  Remark {\ref{H_smaller_1/2}}.

\begin{prop}\label{prop2}
If $H>1/6$ and if $\sigma\in\mathcal{C}^5(\R)$ satisfies a global Lipschitz condition, then the equation
\begin{equation}\label{eq1b=0}
X_t=x_0+\int_0^t \sigma(X_s)d^\circ B_s, \quad t\in [0,1]
\end{equation}
admits a unique
solution $X$ in the set of processes of the form $X_t=f(B_t)$ with $f\in{\mathcal C}^5(\R)$.
The solution is given by 
$
X_t=\phi(x_0,B_t),\,t\in[0,1],
$
where $\phi$ is defined by (\ref{phi}).
\end{prop}
{\bf Proof}. See \cite{lpma1024}, Theorem 2.10. \fin
\end{itemize}

\medskip

\begin{remark}{\label{H_smaller_1/2}}

{\rm 
In \cite{nosi}, Nourdin and Simon developed recently a new concept, namely the 
Newton-C\^otes integral corrected by L\'evy areas, in order to study equation (\ref{eq1})
for any $H\in (0,1)$. It allows to use a fixed point
theorem to obtain  existence and  uniqueness in the set of processes whose paths are H\"older continuous 
of index $\alpha\in(0,1)$ and not only in the more restrictive - and a little arbitrary - set of processes 
of the form $X_t=f(B_t,A_t)$ with $f:\mathbb{R}\times [0,1]\rightarrow\mathbb{R}$ regular
enough and $A$ a process with $\mathcal{C}^1$-trajectories, as shown in \cite{lpma1024}.
}

\end{remark}

\bigskip

Approximation schemes for  stochastic differential equations of the type (\ref{eq1})  are studied
only in few articles, see, e.g., \cite{neuenkirch} and the references therein.
In \cite{cras}, the second-named author considers the approximation of autonomous 
differential equations driven by H\"older continuous
functions (of any fractal index $0<\alpha<1$). He determines upper bounds for the order 
of convergence of the 
Euler scheme and a Milshtein-type scheme, see also \cite{talay}, and applies then his 
results to the case of the fBm. In \cite{neuenkirch},
the first-named author studies the following equation with additive fractional noise
\begin{eqnarray} \label{add_eq}
X_t=x_0+\int_0^t \sigma(s)dB_s +\int_0^t b(s,X_s)ds,\quad t\in [0,1] \end{eqnarray} under the hypothesis $H>1/2$. For a mean-square-$L^{2}$-error criterion, 
he derives by means of the Malliavin calculus
the exact rate of convergence of the Euler scheme, also for non-equidistant discretizations. Moreover, the optimal approximation of equation  (\ref{add_eq}) is also studied in \cite{neuenkirch}.

In this paper, we are interested in the exact
rate of convergence of the Euler scheme associated to (\ref{eq1f})
and of the Crank-Nicholson schemes associated to (\ref{eq1b=0}).
Thus here, compared to \cite{neuenkirch}, we study the non-additive case.
We obtain two types of results (see Section 3 for precise statements):
\begin{description}
\item[(i)]  If $H>1/2$ and under standard assumptions on $\sigma$ and $b$, then the classical Euler scheme
$\overline{X}^n$ with step-size $1/n$ for equation (\ref{eq1f}) defined by
\begin{equation}\label{euler-eq1}
\left\{
\begin{array}{lll}
\overline{X}^n_0=x_0\\
\overline{X}^n_{(k+1)/n}=\overline{X}^n_{k/n}+
\sigma(\overline{X}^n_{k/n})\Delta B_{k/n}
+b(\overline{X}^n_{k/n})\frac{1}{n},\quad 
k\in\{0,\ldots,n-1\},
\end{array}
\right.
\end{equation}
and $\overline{X}^n_t=\overline{X}^n_{[nt]/n}$ for $t\in [0,1]$
verifies
\begin{equation}\label{intro-1}
n^{2H-1}\left[ \, \overline{X}^n_1 - X_1 \, \right] 
\,{\stackrel{{\rm a.s.}}{\longrightarrow}}\, 
-\frac{1}{2}\int_0^1 \sigma'(X_s)D_sX_1ds.
\end{equation}
Here $D_s X_t$, $s,t \in [0,1]$ denotes the Malliavin derivative
at time $s$ of $X_t$ with respect to the fBm $B$. This result is somewhat surprising because
it does not have an analogue in the case of the standard Brownian motion. 
Indeed, in this framework, or more generally when
 SDEs driven by semimartingales are considered, it is generally shown that  $\overline{X}^n_1$ converges a.s. 
to $X_1$ and then that the 
correctly renormalized difference converges in law, see, e.g, \cite{KP} and Remark   {\ref{euler_t1}}, point 3. 
For the approximation of It\^{o}-SDEs with respect to mean square error criterions, see, e.g.,   \cite{Kloeden_Platen}  or \cite{Milstein} and Remark {\ref{euler_t1}}, point 3.
 
Moreover, if we consider  the  global error on the  interval $[0,1]$  of  the Euler scheme, we obtain 
$$
n^{2H-1}\left\| \overline{X}^n - X \right\|_\infty 
\,{\stackrel{{\rm a.s.}}{\longrightarrow}}\, 
\frac{1}{2}\sup_{t\in [0,1]}\left| \int_0^t \sigma'(X_s)D_sX_tds \right|.$$

\item[(ii)]  Assume that $1/6<H<1/2$ and let us consider the  Crank-Nicholson scheme $\widehat{X}^n$ 
with step-size $1/n$ associated to (\ref{eq1b=0}):
\begin{equation}\label{crank-eq1}
\left\{
\begin{array}{lll}
\widehat{X}^n_0=x_0\\
\widehat{X}^n_{(k+1)/n}=\widehat{X}^n_{k/n}+\frac{1}{2}\left(\sigma(\widehat{X}^n_{k/n})
+\sigma(\widehat{X}^n_{(k+1)/n})\right)(B_{(k+1)/n}-B_{k/n}),\\
\hskip10cm k\in\{0,\ldots,n-1\},
\end{array}
\right.
\end{equation}
and $\widehat{X}^n_t=\widehat{X}^n_{[nt]/n}$ for $t\in [0,1]$. Here we obtain the following rates of convergence:
\begin{itemize}
\item {\it (Exact rate)} If the diffusion coefficient $\sigma\in {\mathcal{C}}^1$ 
satisfies 
$$
\sigma(x)^2=\alpha x^2+\beta x+\gamma \mbox{ with some } \alpha, \beta, \gamma\in\R,
$$   we have 
\begin{equation}\label{intro-2}
n^{3H-1/2}\left[ \widehat{X}^n_1 - X_1 \right] 
\,{\stackrel{\mathcal{L}}{\longrightarrow}}\, 
\sigma_H\,\frac{\alpha}{12}\, \sigma(X_1)\, G,
\end{equation}
with $G\sim {\rm N}(0,1)$ independent  of $X_1$ and $\sigma^2_H$ given by (\ref{sigma}).
We prove also an equivalent of (\ref{intro-2}), at the global level:
$$
n^{3H-1/2}\sup_{k\in\{0,\ldots,n\}}\left| \widehat{X}^n_{k/n} - X_{k/n} \right| 
\,{\stackrel{\mathcal{L}}{\longrightarrow}}\, 
\sigma_H\,\frac{\alpha}{12}\,\sup_{t\in [0,1]} \left|\sigma(X_t)\, W_t\right|,
$$
with $W$ a standard Brownian motion independent of $X$.

Compared to the above result for the Euler scheme, 
the convergence  to a mixing law obtained here is classical in the semimartingale
framework. In the fBm framework, such a phenomenon was already obtained in three recent
papers for $1/2<H<3/4$: in \cite{CNW}, the authors study the asymptotic behavior of the power variation
of processes of the form $\int_0^T u_s dB_s$, while, in \cite{LL} and \cite{LL2} the asymptotic
behavior of $\int_0^t f(\overline{X}^n_s)G(\dot{\overline{X}^n_s} \,n^{H-1})ds$ is studied,
where $\overline{X}^n$ denotes the broken-line approximation with stepsize $1/n$ of the solution $X$ of (\ref{eq1b=0}) and $\dot{\overline{X}^n}$ its derivative.

\item {\it (Upper bound)} If $1/3<H<1/2$ and $\sigma \in \mathcal{C}_{b}^{\infty}$  is bounded, we have
\begin{equation}\label{intro-3}
\mbox{for any }\alpha<3H-1/2,\quad n^{\alpha}\left[\widehat{X}^{n}_1-X_1\right]
\,{\stackrel{{\rm Prob}}{\longrightarrow}}\,\,\,
0.
\end{equation}
\end{itemize}
Note that (i) covers in particular the case of a linear diffusion coefficient, 
while in (ii) we consider smooth diffusion coefficients, which are bounded and therefore nonlinear. 
The exact rate of convergence in the general case (ii) remains an open problem although
it seems that it is again $3H-1/2$, as in (i).


\end{description}
The paper is organized as follows. In Section 2, we recall a few facts about the Malliavin calculus
with respect to the fBm $B$. Section 3 contains the results concerning the exact rates of convergence
for the Euler and the Crank-Nicholson schemes associated to (\ref{eq1f}) and 
(\ref{eq1b=0})  respectively. 
The proofs of the results for the Euler scheme are postponed to Section 4. 

\bigskip

\section{Recalls of Malliavin calculus with respect to a fBm}

Let us give a few facts about the Gaussian structure of fBm and its Malliavin derivative process, following Section 3.1 in \cite{nualouk} and Chapter 1.2 in \cite{nual}.
Let $\EE$ be the set of step-functions on $[0,1]$. Consider the Hilbert space $\HH$ defined as the closure of $\EE$ with respect to the scalar product 
$${\lpa\Un_{[0,t]}, \Un_{[0,s]}\rpa}_{\HH}\; =\; R_H(t,s), \quad s,t \in [0,1].$$
More precisely, if we set
$$K_H(t,s)\; =\; {\Gamma\left( H + 1/2\right)}^{-1}{(t-s)}^{H-1/2}F\left( H-1/2,1/2 -H ; H+ 1/2 , 1 - t/s\right),$$
where $F$ denotes  the standard hypergeometric function, 
and if we define the linear operator $K_H^*$ from $\EE$ to $L^2([0,1])$ by
$$(K_H^*\varphi)(s)\; =\; K_H(T,s)\varphi(s) \; + \; \int_s^T(\varphi(r) - \varphi(s))\frac{\partial K_H}{\partial r}(r,s)\, dr, \quad \varphi \in \HH, \,\, s \in [0,1], $$
then $\HH$ is isometric to $L^2([0,1])$ due to the equality
\begin{equation}
\label{HH}
{\lpa\varphi, \rho\rpa}_{\HH}\; = \;\int_0^T  (K_H^*\varphi)(s)(K_H^*\rho)(s)\, ds, \quad \varphi, \rho \in \HH.
\end{equation}
The process $B$ is a centered Gaussian process with covariance function $R_H$, hence its associated Gaussian space is 
isometric to $\HH$ through the mapping $\Un_{[0,t]}\mapsto B_t$. 

\vspace{2mm}

Let $f : \R^n \rightarrow \R$ be a smooth function with compact support and consider the random variable 
$F = f(B_{t_1}, \ldots, B_{t_n})$ (we then say that $F$ is a smooth random variable). The derivative process of $F$ is the 
element of $L^2(\Omega, \HH)$ defined by
$$D_sF\; =\; \sum_{i =1}^n \frac{\partial f}{\partial x_i}(B_{t_1}, \ldots, B_{t_n})
\Un_{[0,t_i]}(s), \quad s \in [0,1].$$
In particular $D_s B_t = \Un_{[0,t]}(s)$. As usual, $\sk^{1,2}$ is 
the closure of the set of smooth random variables with respect to the norm 
$$\lva\lva F\rva\rva_{1,2}^2 \; = \; { \rm E}\lcr |F|^2\rcr + { \rm E} \lcr \lva\lva D.F\rva\rva_{\HH}^2 \rcr.$$
The divergence operator $\delta$ is the adjoint of the derivative operator. If a random variable $u \in L^{2}(\Omega, \HH)$
belongs to the domain of the divergence operator, then  $\delta(u)$ is defined by the duality relationship
$$ { \rm E} (F \delta(u))= { \rm E} \langle D F, u \rangle_{\HH}, $$ for every $F \in \sk^{1,2}$. \\
Finally, let us recall the following result proved in \cite{nosi05}:
\medskip
\begin{prop}\label{prop34}
Let $H>1/2$,  $\sigma\in \mathcal{C}^2_b$ and  $b\in \mathcal{C}^1_b$. Then we have for the  unique solution $X=(X_t,\,t\in[0,1])$ of equation (\ref{eq1f}) that  $X_t\in \sk^{1,2}$ for any $t  \in  (0,1]$ and 
\begin{equation}\label{malliavin}
D_sX_t=\sigma(X_s){\rm exp}\left(\int_s^t b'(X_u)du+\int_s^t \sigma'(X_u)d^-B_u\right), 
\quad 0\leq s \leq t\le 1,\quad t>0.
\end{equation}
\end{prop}

\bigskip

\section{Exact rates of convergence}
In the sequel, we will assume that $[0,1]$ is partitioned by $\{0=t_0<t_1<\ldots<t_n=1\}$ with $t_k=k/n$, $0\le k\le n$.
Rates of convergence will thus be given relative to this partition scheme. For simplicity, we write 
$\Delta B_{k/n}$
instead of $B_{(k+1)/n}-B_{k/n}$.  
For the usage of non-equidistant discretization in the simulation of fBm, see Remark \ref{rm34}, point 2.
\smallskip
\subsection{Euler scheme}
In this section, we assume that $H>1/2$ and we consider equation
(\ref{eq1f}), i.e., the integral with respect to $B$ is defined by (\ref{forward}). 
The Euler scheme $\overline{X}^n$ for equation  (\ref{eq1f}) 
is defined by (\ref{euler-eq1}).
The following theorem shows that the exact rate of convergence of the Euler scheme is in general $n^{-2H+1}$ for the error at the single point $t=1$.

\medskip
\begin{theorem}\label{thm1}
Let  $b\in {\mathcal C}^{2}_b$ and 
$\sigma\in {\mathcal C}^{3}_b$. 
Then, as $n\rightarrow\infty$, we have
\begin{equation}\label{cv1_as}
n^{2H-1}\left[ \overline{X}^n_1 - X_1 \right]  
\,{\stackrel{{\rm a.s.}}{\longrightarrow}}\, 
-\frac{1}{2}\int_0^1 \sigma'(X_s)D_sX_1ds.
\end{equation}
\end{theorem}
\noindent
{\bf Proof}. We postpone it to Section 4. \fin
\noindent
\begin{remark}{\label{euler_t1}}
{\rm \begin{enumerate} 
\item The asymptotic constant of the error does not vanish, e.g., for the linear equation with constant coefficients, i.e.,
\begin{equation}\label{linear}
X_t=x_0+  \gamma \int_0^t X_s d^{-}B_s + \beta \int_0^t X_s ds,\quad t\in [0,1]
\end{equation}
with $\gamma, \beta \in \R$, $\gamma \neq 0$ and $x_{0} \neq 0$. 
\item The appearance of the Malliavin derivative in the asymptotic constant seems to be due to the fact 
that  $(D_{t}X_1)_{t \in [0,1]}$  measures the functional dependency of $X_1$ on the 
driving fBm, see \cite{NuSa}. 
\item  For the  It\^{o}-SDE, i.e., the case $H=1/2$, it is shown in \cite{KP} that
\begin{equation*}  n^{1/2}  \left[ \overline{X}^{n}_{1}-X^{\textrm{It\^{o}}}_{1}  \right ]   \stackrel{\mathcal{L}}{\longrightarrow}  - \frac{1}{\sqrt{2}}Y_{1}  \int_{0}^{1}  \sigma\sigma'(X^{\textrm{It\^{o}}}_{s})Y_{s}^{-1} \, dW_s  \end{equation*}
as $n \rightarrow \infty$, 
with  a Brownian motion  $W$, which is independent of the Brownian motion $B$, and 
$$ Y_{s}= \exp \left(\int_{0}^{s} b'(X^{\textrm{It\^{o}}}_{u}) -\frac{1}{2}\sigma\sigma'(X^{\textrm{It\^{o}}}_{u})\, d u + \int_{0}^{s} \sigma'(X^{\textrm{It\^{o}}}_{u})\, dB_{u}  \right), \quad s \in [0,1].$$ In \cite{Cam_Hu}  the analogous assertion for the mean-square error is established, i.e.,
$$
n  \, { \rm E} \left| \overline{X}^n_1 - X_1^{\textrm{It\^{o}}} \right|^{2} 
\,{\stackrel{}{\longrightarrow}}\, \frac{1}{2}{\rm E} \left| Y_{1} \, \int_{0}^{1}  \sigma\sigma'(X^{\textrm{It\^{o}}}_{s})Y_{s}^{-1} \, dW_s  \right|^{2} $$
 as $n\rightarrow\infty.$

\item 
Assume that  $b\in {\mathcal C}^{2}_b$, 
$\sigma\in {\mathcal C}^{3}_b$  and that  additionally $b$ and $\sigma$ are bounded with  $\inf_{x \in \R} |\sigma(x)| >0$.
Under these
stronger assumptions, which are  a priori  only of technical nature, we can show - by applying the same techniques - that  Theorem {\ref{thm1}}  is also valid with respect to the mean square error, i.e.,
$$
n^{2H-1} \left( { \rm E} \left| \overline{X}^n_1 - X_1 \right|^{2} \right)^{1/2}
\,{\stackrel{}{\longrightarrow}}\, 
\frac{1}{2}\left( { \rm E} \left| \int_0^1 \sigma'(X_s)D_sX_1ds \right|^{2} \right)^{1/2} $$
 as $n\rightarrow\infty$. 
Note that under the above assumptions on $b$ and $\sigma$  equation (\ref{D}) simplifies to $$ 
A'_t=\frac{b(\phi(A_t,B_t))}{\sigma(\phi(A_t,B_{t})) }, \quad A_0=x_0, \qquad t \in [0,1], $$ which  allows us to control the integrability of the remainder terms in the error expansions made in the Proof of Theorem {\ref{thm1}}.

\item For equation (\ref{linear}) and when $H< 3/4$, it is possible to go further and to obtain a
convergence in law for the third term in the asymptotic development of $\overline{X}^n_1$, see also Theorem \ref{thm2}: We have,
as $n\rightarrow +\infty$,
\begin{equation}\label{cvlaw}
\begin{array}{rcl}
\overline{X}^n_1&\,{\stackrel{{\rm a.s.}}{\longrightarrow}}\,&X_1\\
n^{2H-1}\left[\overline{X}^n_1-X_1\right]
&\,{\stackrel{{\rm a.s.}}{\longrightarrow}}\, 
&-\frac{\gamma^2}{2}\,X_1\\
n^{2H-1/2}\left[\overline{X}_1^n-X_1+\frac{\gamma^2}{2}\,X_1\,n^{1-2H}\right]
&\,{\stackrel{\mathcal{L}}{\longrightarrow}}\,&
-\frac{\gamma^2}{2}\, X_1\,G
\end{array}
\end{equation}
with $G$ a centered Gaussian random variable. For $H=3/4$ the last convergence
is again valid if one replaces $n^{2H-1/2}$ by $n^{2H-1/2}\,(\log n)^{-1/2}$. Indeed,
we have 
$$
\overline{X}^n_1=x_0\,{\rm exp}  \left( \gamma B_1+\beta-\frac{1}{2}\gamma^2\,\sum_{k=0}^{n-1}(\Delta B_{k/n})^2+ R_{n}
\right )
$$ with $n^{2H-1/2} |R_{n}| \,  {\stackrel{{\rm a.s.}}{\longrightarrow}} \,  0$ as $n \rightarrow \infty$. 
Thus it holds, as $n\rightarrow +\infty$,
$$
\overline{X}_1^n-X_1+\frac{\gamma^2}{2}\,X_1\,n^{1-2H}\approx -\frac{\gamma^2}{2n^{2H}}\,X_1\,
\sum_{k=0}^{n-1}[(n^H\Delta B_{k/n})^2-1].
$$
Hence Theorem 3 (or Theorem 6 for the case $H=3/4$) in \cite{CNW} allows us to 
obtain the convergence in law in (\ref{cvlaw}).
When $H>3/4$, it seems to be hard to derive a result in law 
since, in this case, arguments used in the proof of Theorem 3 in \cite{CNW} are not valid anymore. Indeed, in this case,
we do not work in a Gaussian framework,
see, e.g.,  Theorem 8 in \cite{CNW}. 

To overcome this problem one can modify the Euler scheme for the linear equation such that the second order quadratic variation of $B$ appears in the error expansion. See, e.g.,  \cite{LL2} for a similar strategy in the case of weighted $p$-variations of fractional diffusions.  
The second order quadratic variation of  fractional Brownian motion is given by $$ V_{n}^{2}(B)=\sum_{k=1}^{n-1}(B_{(k+1)/n} - 2B_{k/n} + B_{(k-1)/n})^2.$$
It is well known, compare for example \cite{BCIJ},  that
\begin{eqnarray*}
 n^{2H-1} V_{n}^{2}(B) & \stackrel{a.s.}{\longrightarrow}&  4-2^{2H} \\ 
 n^{2H-1/2}  V_{n}^{2}(B)- n^{1/2}(4-2^{2H})  & \stackrel{\mathcal{L}}{\longrightarrow}& G,
\end{eqnarray*}
as $n \rightarrow \infty$,
where $G$ is a centered Gaussian random variable with  known variance $c_{H}^{2}>0$. Moreover, by an obvious modification of Proposition \ref{prop3} we also have 
$$ \left(B_{1},  n^{2H-1/2}  V_{n}^{2}(B)- n^{1/2} (4-2^{2H}) \right)   \stackrel{\mathcal{L}}{\longrightarrow}  \left(B_{1}, G \right)$$
as $n \rightarrow \infty$, with $G$ independent of $B_{1}$.
 
For the following approximation scheme for the linear equation (\ref{linear})
\begin{equation*}
\left\{
\begin{array}{lll}
\widetilde{X}^n_0=x_0\\
\widetilde{X}^n_{(k+1)/n}=\widetilde{X}^n_{k/n}+ \gamma \widetilde{X}^n_{k/n}   \Delta_{k/n}B +  \frac{\gamma^{2}}{2}  \widetilde{X}^n_{k/n} \Delta_{k/n}B \, \Delta_{(k-1)/n}B +  \beta\widetilde{X}^n_{k/n}\frac{1}{n},\\
\hskip10cm k\in\{0,\ldots,n-1\},
\end{array} \right.
\end{equation*}
we get by straightforward calculations 
\begin{equation*}
\overline{X}^n_1=X_{1} \,  {\rm exp}  \left(-\frac{1}{4}\gamma^2V_{n}^{2}(B)\, -\frac{\beta^{2}}{2} \frac{1}{n} - \beta \gamma B_{1} \frac{1}{n} + R_{n} \right)
\end{equation*}
with $ n^{ \min \{{2H-1/2, 1} \} } |R_{n}| 
\,  {\stackrel{{\rm a.s.}}{\longrightarrow}} \,  0$ as $n \rightarrow \infty$.
Hence we obtain
$$  n^{2H-1} \left[  \widetilde{X}^n_{1} -X_{1} \right]  \stackrel{a.s.}{\longrightarrow}  - \frac{\gamma^{2}}{4}(4-2^{2H})X_{1}   $$
as $n \rightarrow \infty$.
Furthermore, we get the following error expansions according to the different  values of $H$.
\begin{itemize}
\item[(i)] Case $1/2<H<3/4$:
$$ n^{2H-1/2} \left[\widetilde{X}_1^n-X_1+ \frac{\gamma^{2}}{4}(4-2^{2H})X_{1}n^{1-2H}
 \right ] \stackrel{\mathcal{L}}{\longrightarrow} - \frac{\gamma^{2}}{4} X_{1} G. $$
\item[(ii)] Case $H=3/4$: 
$$ n \left[\widetilde{X}_1^n-X_1+ \frac{\gamma^{2}}{4}(4-2^{3/2})X_{1}n^{-1/2}
 \right ] \stackrel{\mathcal{L}}{\longrightarrow} - \frac{\gamma^{2}}{4} X_{1} G - \frac{\beta^{2}}{2}X_{1} - \beta \gamma X_{1} B_{1}. $$
\item[(iii)] Case $3/4<H<1$:
$$ n \left[\widetilde{X}_1^n-X_1+  \frac{\gamma^{2}}{4}(4-2^{2H})X_{1}n^{1-2H}
 \right ] \stackrel{\mathcal{L}}{\longrightarrow} - \frac{\beta^{2}}{2}X_{1} - \beta \gamma X_{1} B_{1}. $$
\end{itemize}

Thus for this scheme we  get -- according to the values of $H$ -- different error expansions  due to the drift part of the equation. 
\end{enumerate}}
\end{remark}

\bigskip
\bigskip

\noindent
For the global  error  on the interval $[0,1]$, we  obtain the  following result. 
\begin{theorem}\label{thm3}
\begin{equation}\label{cv2_as}
n^{2H-1}   \| \overline{X}^n - X \|_{\infty} 
\,{\stackrel{{\rm a.s.}}{\longrightarrow}}\, 
\frac{1}{2} \sup_{t \in [0,1]} \left| \int_0^t \sigma'(X_s)D_sX_tds \right|.
\end{equation}
\end{theorem}
{\bf Proof}. We postpone it to Section 4. \fin
Hence the Euler scheme obtains   the same exact rate of convergence  for the global error on the interval $[0,1]$  as for the error at the single point $t=1$. 
Moreover we have a.s.
$$ \sup_{t \in [0,1]} \left| \int_{0}^{t} \sigma'(X_{s})D_s X_{1} \, ds \right| = 0 $$ if and only if  a.s.
$$ X_t  \in  (\sigma \sigma')^{-1} (\{ 0 \}) \quad \textrm{for all} \quad t \in [0,1]. $$

\medskip

\medskip

\begin{remark}\label{rm34}
{\rm \begin{enumerate}
\item If  $b=0$,  Theorem \ref{thm1} and \ref{thm3} are again valid under the weaker assumption that
$\sigma\in {\mathcal C}^{1}_b$. Since in this case
$$D_s X_t=\sigma(X_s){\rm exp} \left( \int_s^t \sigma'(X_u)d^-B_u \right)=\sigma(X_t) \qquad s\in [0,t],$$ 
which
is an obvious consequence of the change of variable formula for fBm, we have here
$$  n^{2H-1}\left [  \overline{X}^n_{1} - X_{1} \right ] 
\,{\stackrel{{\rm a.s.}}{\longrightarrow}}\, 
-\frac{1}{2} \sigma(X_1)\,\int_0^1 \sigma'(X_s)ds  $$
and
$$ \, \, \quad n^{2H-1}\left\| \overline{X}^n - X \right\|_\infty 
\,{\stackrel{{\rm a.s.}}{\longrightarrow}}\, 
\frac{1}{2}\sup_{t\in [0,1]}\left|\sigma(X_t)\,\int_0^t \sigma'(X_s)ds\right|,$$
respectively.

\item For $H \neq 1/2$ the increments of fractional Brownian motion are  correlated. Therefore  the exact simulation of $B(t_{1}), \ldots, B(t_{n})$ is in general computationally very expensive. The Cholesky decomposition method, which is to our best knowledge the only known exact method for the non-equidistant simulation of fractional Brownian,   requires $O(n^{3})$ operations. Moreover the covariance matrix, which has to be decomposed,  is ill-conditioned.
If the discretization is equidistant, i.e., $t_{i}=i/n$, $i=1, \ldots, n$,
 the computational cost can be lowered considerably, making use of the stationarity of the increments of fractional Brownian motion. For example, the Davies-Harte algorithm for the equidistant simulation of fractional Brownian motion has 
computational cost $O(n \log(n))$, see, e.g., \cite{Craig}. For a comprehensive survey of simulation methods for fractional Brownian motion, see, e.g., \cite{Coeur}.
\item In the  Skorohod setting, it is in general  difficult to write an
Euler type scheme associated to  equation (\ref{eq1}), even if $H> 1/2$ and $b=0$. 
Indeed, in this case, by using the 
integration by parts rule $\delta(Fu)=F\delta(u)-\langle DF,u
\rangle_{\cal{H}}$ for the Skorohod integral  and by approximating
$X_{(k+1)/n}$ by $X_{k/n}+\int_{k/n}^{(k+1)/n} \sigma(X_{k/n})\delta
B_s$ (as in the case $H=1/2$), one obtains
$$
\overline{X}^{n}_{(k+1)/n}=\overline{X}^{n}_{k/n}
+\sigma( \overline{X}^{n}_{k/n}) \left(B_{(k+1)/n}-B_{k/n}\right) - \sigma'(\overline{X}^{n}_{k/n})
\langle D\overline{X}^{n}_{k/n},1_{[k/n,(k+1)/n]}\rangle_{\mathcal{H}}.
$$ 
The problem is that the Malliavin derivative $D\overline{X}^{n}_{k/n}$ appears, which is  difficult to compute  directly. Moreover, the error analysis of such an approximation seems also to be very difficult,  because the $L^{2}$-norm of the Skorohod
integral involves the first Malliavin derivative of the integrand. Thus, for analyzing such an approximation scheme, we need also to control the difference between the Malliavin derivative of the solution and  the Malliavin derivative of the approximation. But this  involves the
second Malliavin derivative etc. and we cannot have closable
formulas. It is one of the reasons for which we preferred here to work within the Russo-Vallois framework, instead of
the Skorohod one. Another reason is that the Russo-Vallois framework is,
from our point of view, simpler in the one-dimensional  case than the Skorohod one, as it is shown 
in \cite{lpma1024}.
\end{enumerate}}
\end{remark}

\bigskip

\bigskip

\subsection{Crank-Nicholson scheme}
In this section, we assume that $1/6<H<1/2$ and we consider equation (\ref{eq1b=0}). 
Let $\widehat{X}^n$ be the
Crank-Nicholson scheme defined by (\ref{crank-eq1}), which is
the canonical scheme associated 
to equation (\ref{eq1b=0}), since the integral with respect to the driving fBm $B$ is 
defined by the symmetric integral. 
It is an implicit scheme, but it is nevertheless well-defined, since
for $n$ sufficiently large, $x\mapsto x-\frac{1}{2}\Delta B_{k/n}\,\sigma(x)$ is 
invertible.
Although (\ref{crank-eq1}) seems to be rather close to (\ref{euler-eq1}) with $b= 0$, 
the situation is in fact here significantly more difficult.  
That is why we study  the rate of convergence for the Crank-Nicholson
scheme only in the following particular cases: 
\begin{itemize}
\item Case 1: $1/6<H<1/2$ and $\sigma \in \mathcal{C}^{1}$ satisfies
$\sigma(x)^2 = \alpha x^2+ \beta x+ \gamma$ for some $\alpha,\beta, \gamma\in\R$,\item Case 2: $1/3<H<1/2$ and $\sigma \in \mathcal{C}^{\infty}_{b}$   bounded. 
\end{itemize}

\subsubsection{Case 1} Compared to 
Theorem \ref{thm1}, 
we have here a convergence in law. Moreover, the limit of the error  is expressed as a mixed law between 
$B_1$ and an 
independent standard Gaussian 
random variable $G$,  see also Remark \ref{euler_t1}, point 4. 

\medskip
\begin{theorem}\label{thm2}
Assume that $1/6<H<1/2$ and  $\sigma\in {\mathcal C}^1$ satisfies $\sigma(x)^2= \alpha x^2+ \beta x+ \gamma$ for some $\alpha,\beta,\gamma\in\R$.
Then, as $n\rightarrow\infty$, we have
\begin{equation}\label{cv2}
n^{3H-1/2}\left[ \widehat{X}^n_{1} - X_{1} \right] 
\,{\stackrel{\mathcal{L}}{\longrightarrow}}\, 
\sigma_H\,\frac{\alpha}{12} \,\sigma(X_1)\,G.
\end{equation}
Here $G\sim {\rm N}(0,1)$ is independent  of $X_1$ and
\begin{equation}\label{sigma}
\sigma_H^2=4/3+1/3\sum_{\ell=1}^\infty
\theta(\ell)^3,
\mbox{ where }2\theta(\ell)=(\ell+1)^{2H}+(\ell-1)^{2H}-2\ell^{2H}.
\end{equation}
\end{theorem}

In fact, we have also a result at a functional level:
\begin{theorem}\label{thm4}
Assume that $1/6<H<1/2$ and  $\sigma \in \mathcal{C}^{1}$ satisfies $\sigma(x)^2=\alpha x^2+ \beta x+ \gamma$ for some $\alpha,\beta,\gamma\in\R$.
Then, as $n\rightarrow\infty$, we have
\begin{equation}\label{cv4}
n^{3H-1/2}\sup_{k\in\{0,\ldots,n\}}\left|\widehat{X}^n_{k/n} - X_{k/n} \right| 
\,{\stackrel{\mathcal{L}}{\longrightarrow}}\, 
\sigma_H\,\frac{\alpha}{12}\,\sup_{t\in [0,1]} \left|\sigma(X_t)\, W_t\right|.
\end{equation}
Here $W$ is a standard Brownian motion independent of $X$ and $\sigma_H>0$ is once again 
given by (\ref{sigma}).
\end{theorem}

\medskip

\begin{remark}\label{rk2}
{\rm 
\begin{enumerate}
\item In \cite{lpma1024}, equation (\ref{eq1b=0}) is also studied and  
it is shown that $\widehat{X}^{n}_{1}$ converges in probability if and only if $H>1/6$ and that, 
in this case, the limit is $X_{1}$, the solution of equation  (\ref{eq1b=0}) at $t=1$.
Of course, this fact is also an obvious consequence of Theorem \ref{thm2}.
\item Let us show how the constant in (\ref{sigma}) appears. Set 
$Y_n=n^{-1/2}\sum_{k=0}^{n-1} (n^H\,\Delta B_{k/n})^{3}$. We claim that ${\rm Var}[Y_n]\rightarrow \sigma_H^2$
as $n\rightarrow +\infty$. Indeed, using $(n^HB_{t/n})_{t \in [0, \infty)}\,{\stackrel{\mathcal{L}}{=}}\,(B_t)_{t \in [0, \infty)}$, we have 
$$
{\rm Var}[Y_n]=n^{-1}\sum_{k,\ell=0}^{n-1} { \rm E} \left[(B_{k+1}-B_k)^3(B_{\ell+1}-B_\ell)^3\right].
$$
Since $x^3=H_3(x)+3H_1(x)$, where
$H_1(x)=x$ and $H_3(x)=x^3-3x$
denote the first and third Hermite polynomial, we can write, by using the well-known identity ${ \rm E} [H_i(X)H_j(Y)]=0$ for $i \neq j$ 
and
${ \rm E}  [H_i(X)H_i(Y)]= { \rm E} [XY]^i/i!$  for a centered Gaussian vector  $(X,Y)$ with 
${ \rm Var}[X]={ \rm Var}[Y]=1$ (see, e.g., \cite{nual}, Lemma 1.1.1): 
$$
{\rm Var}[Y_n]=(6n)^{-1}\sum_{k,\ell=0}^{n-1}\theta(\ell-k)^3+9n^{-1}\sum_{k,\ell=0}^{n-1}\theta(\ell-k), 
$$
with 
$$
2\theta(\ell-k)=2 {\rm E} \left[(B_{k+1}-B_k)(B_{\ell+1}-B_\ell)\right]=
|\ell-k+1|^{2H}+|\ell-k-1|^{2H}-2|\ell-k|^{2H}.
$$
On the one hand, let us remark that $n^{-1}\sum_{k,\ell=0}^{n-1}\theta(\ell-k)=n^{-1}
 {\rm E} \left(\sum_{k=0}^{n-1} B_{k+1}-B_k \right)^{2}=n^{2H-1}\longrightarrow 0$ as $n\rightarrow +\infty$, 
if $H<1/2$.  
On the other hand, we can write 
$$
\sum_{k,\ell=0}^{n-1}\theta(\ell-k)^3=\sum_{k=0}^{n-1}\theta(0)^3+2\sum_{k=0}^{n-1}\sum_{\ell=k+1}^{n-1}\theta(\ell-k)^3
=8n+2\sum_{k=0}^{n-1}\sum_{\ell=1}^{n-k-1}\theta(\ell)^3.
$$
Consequently, since $\theta(\ell) < 0$ for $H<1/2$,  we deduce by using Cesaro's theorem that $ { \rm Var}[Y_n]\longrightarrow \sigma_H^2$
given by (\ref{sigma}), as $n\rightarrow +\infty$.
\item   In \cite{huesler},  in particular the approximation of fBm by its piecewise linear interpolation 
$$ \widetilde{B}^{n}_t = B_{k/n} + (nt-[nt]) (B_{(k+1)/n} -B_{k/n}), \qquad t \in [k/n, (k+1)/n]$$
 is studied. It is  shown that the correct renormalization of  $\| B- \widetilde{B}^{n}\|_{\infty}$    converges to the Gumbel distribution, i.e., 
$$ P \left( \| \widetilde{B}^{n} -B \|_{\infty} \leq \sigma_{n}(\nu_{n} + x/ \nu_{n}) \right) \longrightarrow \exp(-\exp(-x))$$
as $n \rightarrow$ for $x \in \R$, where $ \sigma_{n} \approx c_{H} n^{-H}$ with $c_{H}>0$  and $\nu_{n}$ is in terms of $\log(n)$. Since $3H-1/2 > H$ for $H>1/4$, the  analogue of Theorem 2 in the setting of Theorem 4, i.e., 
$$ n^{3H-1/2}\|\widehat{X}^n - X \|_{\infty} 
\,{\stackrel{\mathcal{L}}{\longrightarrow}}\,\sigma_H\,\frac{\alpha}{12}\,\sup_{t\in [0,1]} \left|\sigma(X_t)\, W_t\right| , $$ as $n \rightarrow \infty$,
can  not hold without further restriction of the Hurst parameter.

\end{enumerate}
}
\end{remark}

\medskip

\noindent
For the proof of Theorem \ref{thm2}, we need the following Lemma.
\begin{lemma}\label{lm}
\begin{description}
\item {\it i)} We have, for $H<1/2$,
\begin{equation}\label{lm-1}
\left(B_1,n^{3H-1/2}\sum_{k=0}^{n-1} (\Delta B_{k/n})^{3}\right)\,{\stackrel{\mathcal{L}}{\longrightarrow}}\,  (B_1,G),\mbox{ as }
n\rightarrow +\infty,
\end{equation}
where $G$ is a centered Gaussian random variable with variance $\sigma^2_H$ given by (\ref{sigma}),
independent of $B_1$. 
\item {\it ii)} We have
$$
n^{5H-1/2}\sum_{k=0}^{n-1} (\Delta B_{k/n})^{5}\,{\stackrel{\mathcal{L}}{\longrightarrow}}\,  G'\mbox{, as }n\rightarrow\infty,
$$
where $G'$ is a centered Gaussian random variable.
\item {\it iii)} We have
$$
n^{6H-1}\sum_{k=0}^{n-1} (\Delta B_{k/n})^{6}\,\,\,\,{\stackrel{{\rm Prob}}{\longrightarrow}}\,\,\,\, 15,\mbox{ as }
n\rightarrow +\infty.
$$
\end{description}
\end{lemma}
{\bf Proof of Lemma \ref{lm}}. The second and the third point are classical: we refer to
\cite{BM}. Thus we have only to prove the first point. 
Let us denote by $H_1(x)=x$ and $H_3(x)=x^3-3x$
the first and third Hermite polynomial.
Since $H<1/2$, we have 
$$
n^{-1/2}\sum_{k=0}^{n-1}H_1(n^{H}\Delta B_{k/n})
=n^{H-1/2}B_1 \, \stackrel{\textrm{Prob}}{\longrightarrow} \,  0
$$
and we deduce that the convergence in law (\ref{lm-1}) will hold if and only if
\begin{equation}\label{lm-1ter}
\left(B_1,n^{-1/2}\sum_{k=0}^{n-1} H_3(n^{H}\Delta B_{k/n})\right)
\,{\stackrel{\mathcal{L}}{\longrightarrow}}\, 
(B_1,G),
\end{equation}
as
$n\rightarrow +\infty.$
In \cite{BM}, it is shown that $n^{-1/2}\sum_{k=0}^{n-1} H_3(n^{H}\Delta B_{k/n})
\,{\stackrel{\mathcal{L}}{\longrightarrow}}\,
G$. See also the second point of Remark \ref{rk2}. Then, the proof (\ref{lm-1ter}) is finished by the following Proposition, which is  an obvious consequence
of the main result contained in \cite{PT}.\fin

\begin{prop}\label{prop3}
Let $(F_1^n,F_3^n)$ be a random vector such that, for every $n$, $F_i^n$ ($i=1,3$) is in the $i$-th Wiener
chaos associated to the fBm $B$. If $F_1^n\,{\stackrel{\mathcal{L}}{\longrightarrow}}\,G_1$ and 
$F_3^n\,{\stackrel{\mathcal{L}}{\longrightarrow}}\,G_3$ with $G_i$ ($i=1,3$) some Gaussian variables, then
$\mathcal{L}(F_1^n,F_3^n)\longrightarrow \mathcal{L}(G_1)\otimes \mathcal{L}(G_3)$.
\end{prop}
\medskip
\noindent
{\bf Proof of Theorem \ref{thm2}}. In \cite{lpma1024}, the second-named author proved
\begin{equation}\label{xn1}
\widehat{X}^n_1 = \phi\left( B_1 + \frac{\alpha}{12} \sum_{k=0}^{n-1} (\Delta B_{k/n})^3
+ \frac{\alpha^2}{80} \sum_{k=0}^{n-1}(\Delta B_{k/n})^5 + O\left(\sum_{k=0}^{n-1} (\Delta B_{k/n})^6\right),x_0\right).
\end{equation} 
Note that  $X_{1}=\phi(x_{0},B_{1})$ and $\frac{ \partial \phi}{ \partial x_{2}}(x_{1},x_{s}) = 
\sigma( \phi(x_{1},x_{2}))$.
Consequently by a  Taylor expansion $\widehat{X}^n_1  -X_1$ equals
$$
 \frac{\alpha}{12}\sigma(X_1)\sum_{k=0}^{n-1} (\Delta B_{k/n})^3
+\frac{\alpha^2}{80}\sigma(X_1)\sum_{k=0}^{n-1} (\Delta B_{k/n})^5 +  O \left( |X_{1}| \times \sum_{k=0}^{n-1}(\Delta B_{k/n})^6  \wedge \left(\sum_{k=0}^{n-1}(\Delta B_{k/n})^3 \right)^{2} \right).$$
By the second point of Lemma \ref{lm}, we deduce that
$n^{3H-1/2}\sum_{k=0}^{n-1} (\Delta B_{k/n})^5\longrightarrow 0$ in law, hence in probability.
By the third point of  Lemma \ref{lm}, and since $H>1/6$, we have that
$n^{3H-1/2}\sum_{k=0}^{n-1} (\Delta B_{k/n})^6
=n^{1/2-3H}\times n^{6H-1}\sum_{k=0}^{n-1} (\Delta B_{k/n})^6\longrightarrow 0$ in probability.
By the first point of  Lemma \ref{lm}, and since $H>1/6$, we have that
$n^{3H-1/2}\left(\sum_{k=0}^{n-1} (\Delta B_{k/n})^3\right)^2
=n^{1/2-3H}\left(n^{3H-1/2}\sum_{k=0}^{n-1} (\Delta B_{k/n})^3\right)^2\longrightarrow 0$ in probability.
Then, using again the first point of  Lemma \ref{lm} and Slutsky's Lemma, we obtain (\ref{cv2}).\fin
\medskip
\noindent
{\bf Proof of Theorem \ref{thm4}}. Exactly as in the proof of Theorem 3 in \cite{CNW}, we can prove
$$
\left(B_t,n^{3H-1/2}V_3^n(B)_t\right)\,{\stackrel{\mathcal{L}}{\longrightarrow}}\,  (B_t,\sigma_H\,W_t),\mbox{ as }
n\rightarrow +\infty,$$
in the space $\mathcal{D}([0,1])^2$ equipped with the Skorohod topology. Here, $W$ is a standard
Brownian motion independent of $B$ and 
$V_3^n(B)$ is defined by $V_3^n(B)_{t}=\sum_{\ell=0}^{[nt]-1} (\Delta B_{\ell/n})^{3}$.
To obtain Theorem \ref{thm4}, it suffices then to adapt the proof of Theorem \ref{thm2} below.
\fin

\subsubsection{Case 2}
Now assume that $\sigma \in \mathcal{C}_{b}^{\infty}$ is bounded and that $1/3<H<1/2$. 

In the sequel, we will need some fine properties concerning the $m$-order variation of $B$ 
on the interval $[0,1]$. Let us state them in the following proposition:
\begin{prop}\label{prop1}
Let $h \in \mathcal{C}^{1}_{b}$.
\begin{enumerate}
\item If $m\in\N$ is even then, for any $H\in (0,1)$:
$$
n^{mH-1}\sum_{k=0}^{n-1} h(B_{k/n})(B_{(k+1)/n}-B_{k/n})^{m}
\,\,\,{\stackrel{{\rm Prob}}{\longrightarrow}}\,\,\,
\frac{m!}{2^{m/2}(m/2)!}\int_0^1 h(B_s) ds,\mbox{ as }n\rightarrow+\infty.
$$
\item If $m\in\N\setminus\{1\}$ is odd  then, for any $H\in (1/4,1/2)$ and $\alpha<mH-1/2$:
\begin{equation}\label{cvlw}
n^{\alpha}\,\sum_{k=0}^{n-1} 
\left[h(B_{k/n})+\frac{1}{2}\,h'(B_{k/n})(B_{(k+1)/n}-B_{k/n})\right]\,(B_{(k+1)/n}-B_{k/n})^{m}
\,\,\,{\stackrel{{\rm Prob}}{\longrightarrow}}\,\,\,
0,
\end{equation}
as $n\rightarrow+\infty$.
\end{enumerate}
\end{prop}
{\bf Proof.}
\begin{enumerate}
\item When $h\equiv 1$, it is a classical result: we refer to \cite{BM} for instance. 
To obtain the general case, it suffices to adapt the methodology developed in step 5 of
\cite{GN}, p. 8. or in the proof of Theorem 1 in \cite{CNW}.   
\item Using the same linear regression as in the proof of Theorem 4.1 in  \cite{GNRV}, we can prove that, when $H<1/2$:
$$
n^{\alpha}\,\sum_{k=0}^{n-1} 
[h(B_{(k+1)/n})+h(B_{k/n})]\,(B_{(k+1)/n}-B_{k/n})^{m}
$$
converges in probability to 0, for any $\alpha<mH-1/2$.
Convergence (\ref{cvlw}) can then be obtained using a Taylor expansion and the fact that $H>1/4$.
\end{enumerate}
Details are left to the reader.
\fin

Using the above  proposition we can show the following result.
\begin{theorem}\label{main-result}
Assume that $H\in(1/3,1/2)$ and $\sigma \in \mathcal{C}^{\infty}_{b}$ is bounded. 
Then we have:
\begin{equation}\label{thm-eq2}
\mbox{For any }\alpha<3H-1/2,\quad n^{\alpha}\,\left[\widehat{X}^{(n)}_1-X_1\right]
\,{\stackrel{{\rm Prob}}{\longrightarrow}}\,\,\,
0
\end{equation}
as $n \rightarrow \infty$.
\end{theorem}
%
%
%
{ \bf Proof of Theorem \ref{main-result}.}
In the following, denote 
$$\Delta^j Z_{k/n}=(Z_{(k+1)/n}-Z_{k/n})^j$$
for $j,n\in\N$, $k\in\{0,\ldots,n-1\}$ and a process $Z=(Z_t)_{t\in[0,1]}$. 
When $j=1$, we prefer the notation $\Delta Z_{k/n}$ instead of $\Delta^1 Z_{k/n}$
for simplicity. 
Denote also, for $p \in \N$, 
$$ \Delta^{p}(B) = \max_{k=0, \ldots, n-1} |\Delta^{p}B_{k/n}|.$$ 
Moreover, it is recall  that $\phi$ given by (\ref{phi}) verifies the semigroup property: 
\begin{equation}\label{sgp}
\forall x,y,z\in \R,\quad
\phi(\phi(x,y),z)=\phi(x,y+z)
\end{equation}
and we have \begin{eqnarray*} X_{t}=\phi(x,B_{t}), \qquad t \in [0,1]. \end{eqnarray*}
Simple but tedious computations (see, for instance, \cite{lpma1024}) allow us to obtain:
\begin{equation*}
{\widehat X}^{(n)}_{(k+1)/n}=
\phi\left({\widehat X}^{(n)}_{k/n},
\Delta B_{k/n}
+
f({\widehat X}^{(n)}_{k/n})\Delta^3 B_{k/n}
+g
({\widehat X}^{(n)}_{k/n})\Delta^4 B_{k/n} + O(\Delta^{5}(B)) \right),
\end{equation*}
with $f=\frac{(\sigma^2)''}{24}$ and $g=\frac{\sigma(\sigma^2)'''}{12}$.
We deduce from the semi-group property (\ref{sgp}) that
for every $\ell\in \{1,\ldots,n\}$:
\begin{equation}\label{use-flow}
{\widehat X}^{(n)}_{\ell/n}=
\phi\left(x,
B_{\ell/n}
+\sum_{k=0}^{\ell-1}f({\widehat X}^{(n)}_{k/n})\Delta^3 B_{k/n}
+\sum_{k=0}^{\ell-1}g({\widehat X}^{(n)}_{k/n})\Delta^4 B_{k/n} +  O(n \Delta^{5} (B)   )
\right). 
\end{equation}
In particular, we have
\begin{equation}\label{wearat}
\sup_{\ell\in\{0,\ldots,n\}}\left|\widehat{X}_{\ell/n}^{(n)}-X_{\ell/n}\right|
=\sup_{\ell\in\{0,\ldots,n\}}\left|\widehat{X}_{\ell/n}^{(n)}-\phi(x,B_{\ell/n})\right| = O ( n \Delta^{3} (B)  )
\end{equation}
and then (\ref{use-flow}) becomes
\begin{equation}\label{use-flow2}
{\widehat X}^{(n)}_{\ell/n}=
\phi\left(x,
B_{\ell/n}
+\sum_{k=0}^{\ell-1}f({\widehat X}^{(n)}_{k/n})\Delta^3 B_{k/n}
+\sum_{k=0}^{\ell-1}g(X_{k/n})\Delta^4 B_{k/n} + O ( n^{2} \Delta^{7} (B)  )
\right).
\end{equation}
Due to the assumptions on $\sigma$ we have  that
\begin{equation*}
 \phi(x,y_{2})= \phi(x,y_{1}) +  \sum_{j=1}^{m}\frac{1}{j!} \frac{ \partial^{j} \phi}{ (\partial x_{2})^{j}} (x,y_{1}) (y_{2}-y_{1})^{j} + O((y_{2}-y_{1})^{m+1}  ). \end{equation*}
Thus  we get 
\begin{eqnarray}{\label{eq_re_in_1}}
&& {\widehat X}^{(n)}_{k/n}=
X_{k/n} + \sum_{j=1}^{m} 
\frac{1}{j!} \frac{\partial^{j} \phi}{(\partial x_{2})^{j}}(x, B_{k/n}) \\ && \quad  \times
 \left( \sum_{k_1=0}^{k-1}f(\widehat{X}^{(n)}_{k_1/n})
\Delta^3 B_{k_1/n}  + \sum_{k_1=0}^{k-1}g(X_{k_1/n})
\Delta^4 B_{k_1/n} +O ( n^{2} \Delta^{7} (B)  )
\right)^{j} \nonumber  \\ && \qquad \qquad \qquad + O(n^{m+1} \Delta^{3(m+1)}(B)). \nonumber
\end{eqnarray} 
{(i)} Now assume  for a moment that $H>5/12$. By using (\ref{eq_re_in_1}) with $m=1$ 
and $\frac{\partial \phi }{\partial x_{2}} (x_{1},x_{2}) = \sigma(\phi(x_{1},x_{2})) $  
we get
\begin{eqnarray*} 
{\widehat X}^{(n)}_{k/n}=
X_{k/n}
+\sigma(X_{k/n})\sum_{k_1=0}^{k-1}f(\widehat{X}^{n}_{k_1/n})
\Delta^3 B_{k_1/n}+ O(n^{2} \Delta^{6}(B)).
\end{eqnarray*}
and
\begin{eqnarray} {\label{eq_re_in_2}}
{\widehat X}^{(n)}_{k/n}=
X_{k/n}
+\sigma(X_{k/n})\sum_{k_1=0}^{k-1}f(X_{k_1/n})
\Delta^3 B_{k_1/n}+ O(n^{2} \Delta^{6}(B)).
\end{eqnarray}
 By inserting the previous equality in (\ref{use-flow2}) with $\ell=n$, we obtain
\begin{eqnarray}\label{eq34}  \nonumber  &&  \widehat{X}^{(n)}_{1}=
\phi\left(x,B_{1} +\sum_{k=0}^{n-1}f(X_{k/n})\Delta^3 B_{k/n} + \sum_{k=0}^{n-1}f'\sigma(X_{k/n}) \Delta^3 B_{k/n}\sum_{k_1=0}^k f(X_{k_1/n})\Delta^3 B_{k_1/n}
 \right. \\  && \qquad \qquad   \qquad  \quad \,\,
  \left. +\sum_{k=0}^{n-1}g(X_{k/n})\Delta^4 B_{k/n} + O(n^{3}\Delta^{9}(B))
\right).
\end{eqnarray}
But we have, due to the second point of Theorem \ref{prop1} and the
fact that $g=2\,\sigma f'$ and $X_t=\phi(x,B_t)$:
$$
n^{\alpha}\left(\sum_{k=0}^{n-1}f(X_{k/n})\Delta^3 B_{k/n}
+\sum_{k=0}^{n-1}g(X_{k/n})\Delta^4 B_{k/n}\right)
\,\,\,{\stackrel{{\rm Prob}}{\longrightarrow}}
\,\,\,
0.
$$ 
Moreover, since $H>5/12$ and $\alpha<3H-1/2$,  we have
$$ n^{\alpha+3}\Delta^{9}(B)   \stackrel{{\rm a.s. }}{\longrightarrow} 0. $$
On the other hand, since 
\begin{eqnarray*}
&& 2\sum_{k=0}^{n-1}f'\sigma(X_{k/n})\Delta^3 B_{k/n}\sum_{k_1=0}^k f(X_{k_1/n})\Delta^3 B_{k_1/n}
\\ && \qquad =
\left(\sum_{k=0}^{n-1}f'\sigma(X_{k/n})
\Delta^3 B_{k/n}\right)
\left(\sum_{k_1=0}^{n-1} f(X_{k_1/n})\Delta^3 B_{k_1/n}\right)
-
\sum_{k=0}^{n-1}ff'\sigma(X_{k/n})\Delta^6 B_{k/n}
\end{eqnarray*}
we deduce, this time due to the first and the second point of Theorem \ref{prop1}, that
$$
n^{\alpha}\sum_{k=0}^{n-1}f'\sigma(X_{k/n})
\Delta^3 B_{k/n}\sum_{k_1=0}^k f(X_{k_1/n})\Delta^3 B_{k_1/n}
\,\,\,{\stackrel{{\rm Prob}}{\longrightarrow}}\,\,\,0.
$$
Finally,  we obtain (\ref{thm-eq2})
when $H>5/12$. \\
{ (ii)} To prove the announced result, that is (\ref{thm-eq2}) for arbitrary $H\in(1/3,1/2)$, 
it suffices to use (\ref{eq_re_in_1}) with the appropriate $m$ for the considered $H$ and then to proceed as in {(i)}. 
The remaining details are left to the reader. \begin{flushright} $\square$ \end{flushright}

\section{Proof of Theorems \ref{thm1} and \ref{thm3}}
Throughout this section we assume that  $b\in {\mathcal C}^{2}_b$,  
$\sigma\in {\mathcal C}^{3}_b$ and $H>1/2$.
For  $g:[0,1] \rightarrow \R$ and $\lambda \in (0,1)$ we will use the usual notations
$$
\| g \|_{\infty} = \sup_{t \in [0,1] }|g(t)|, 
\qquad  \|g\|_{\lambda} = \sup_{s,t \in [0,1], s \neq t} \frac{|g(t)-g(s)|}{|t-s|^{\lambda}}.
$$
Moreover positive constants,  depending only on $b$, $\sigma$, their derivatives, $x_{0}$ and $H$,
will be denoted by $c$, regardless of their value. We will write $\Delta$ instead of $1/n$.\\ 

\noindent
The following properties of the function $\phi$ are taken from Lemma 2.1 in \cite{talay}.
\medskip
\begin{lemma}{\label{prop_h}}  Let $\phi$ given by (\ref{phi}). Then we have 
\begin{eqnarray*}
 \textrm{(a)} && \quad \phi(x_{1},x_{2})=\phi(\phi(x_{1},y),x_{2}-y), \qquad x_{1}=\phi(\phi(x_{1},x_{2}),-x_{2}) \\
 \textrm{(b)} && \quad \frac{ \partial \phi}{\partial x_{2}}(x_{1},x_{2}) = \sigma(x_{1})  \frac{ \partial \phi}{\partial x_{1}}(x_{1},x_{2}) \\
\textrm{(c)} && \quad  \sigma^{2}(x_{1})  \frac{ \partial^{2} \phi}{\partial x_{1}^{2}}(x_{1},-x_{2}) -2 \sigma(x_{1}) \frac{\partial^{2} \phi}{\partial x_{1} \partial x_{2}}(x_{1},-x_{2})  + (\sigma \sigma')(x_{1})  \frac{ \partial \phi}{\partial x_{1}} (x_{1},-x_{2}) \\ && \qquad \qquad  \qquad \qquad+ \frac{ \partial^{2} \phi}{\partial x_{2}^{2}}(x_{1},-x_{2})= 0 \\
 \textrm{(d)} &&  \quad 1= \frac{ \partial \phi}{\partial x_{1}}(\phi(x_{1},x_{2}),-x_{2})   \frac{ \partial \phi}{\partial x_{1}}(x_{1},x_{2}) \\ 
\textrm{(e)} && \quad \frac{ \partial \phi}{\partial x_{1}}(x_{1},x_{2}) = \exp \left( \int_{0}^{x_{2}} \sigma'(\phi(x_{1},s)) \, ds \right)
\end{eqnarray*}
for all  $x_{1},x_{2},y \in \R$.
\end{lemma}
\bigskip
The following is  well known and easy to prove.
\medskip
\begin{lemma}\label{rek} Let $n \in \N$ and $a_{i}, b_{i} \in \R$ for $i=1, \ldots, n$. \\ 
(a) For $x_{j}$, $j=0, \ldots, n$ given by the recursion
$$
 x_{j+1} = x_{j}   b_{j} + a_{j}, \quad \quad \quad \quad j=0, \ldots, n-1,
$$  with $x_{0} = 0$,
we have
$$ x_{j} = \sum_{i=0}^{j-1} \, \, a_{i}   \prod_{k=i+1}^{j-1} b_k , \qquad j=1, \ldots n. $$
(b) If
$$
  |x_{j+1}|  \leq |x_{j}|   |b_{j}| + |a_{j}|, \quad \quad \quad \quad j=0,
\ldots, n-1, $$  with $x_{0} = 0$ and $|b_j|\ge 1$ for all $j=0,\ldots,n-1$, 
then
\[ \max_{j=0, \ldots, n} |x_{j}| \, \leq \,  \sum_{i=0}^{n-1} \, \, |a_{i}|    \prod_{k=1}^{n-1} |b_k|.  \]

\end{lemma}
\bigskip
We  will also require that the Euler approximation of the solution and the process $(A_{t})_{t \in [0,1]}$ given by (\ref{D}) can be uniformly bounded in terms of the driving fBm. 

\medskip
\begin{lemma} \label{bound_euler}
We have 
\begin{eqnarray}
 && \sup_{n \in \N} \sup_{k=0, \ldots, n} | \XE_{k/n} |    \, \leq  \,  \exp(c \,  \exp( c\,  (\| B \|_{\infty} + \| B \|_{1/2}^{2} ))) \qquad \textrm{a.s.}, \\ 
 && \qquad \qquad \quad \|A \|_{\infty}  \, \leq \,    \exp(c  \,  \exp( c \,   \| B \|_{ \infty}))  \quad  \qquad \qquad \qquad \textrm{a.s.} \end{eqnarray} \end{lemma}
{\bf Proof.}
We prove only the first assertion, following the proof of Lemma 2.4. in \cite {talay}. The second assertion can be  obtained by a straightforward  application of  Gronwall's Lemma to  equation (\ref{D}).
By Lemma \ref{prop_h} (a)  we have
$$ \phi( \XE_{k/n},-B_{k/n})=\phi(\phi( \XE_{k/n},\Delta B_{k/n}),-B_{(k+1)/n}).$$
Using this, we obtain by the mean value theorem 
$$  \phi( \XE_{(k+1)/n},-B_{(k+1)/n})-  \phi( \XE_{k/n},-B_{k/n}) =  \left[  \XE_{(k+1)/n} - \phi( \XE_{k/n},\Delta B_{k/n} ) \right]   \frac{\partial \phi}{\partial x_{1}} (\xi_{k}, -B_{(k+1)/n})  $$
with $\xi_{k}$ between $ \XE_{(k+1)/n}$ and $\phi( \XE_{k/n}, \Delta B_{k/n})$.
Moreover
$$ \phi( \XE_{k/n},  \Delta B_{k/n} ) = \phi( \XE_{k/n}, 0 ) +  \Delta B_{k/n}    \frac{\partial \phi}{\partial x_{2}} ( \XE_{k/n}, 0)    + \frac{1}{2}    (\Delta B_{k/n} )^{2}    \frac{\partial^{2} \phi}{\partial x_{2}^{2}} ( \XE_{k/n}, \zeta_{k}) $$ 
with $|\zeta_{k}| \leq |\Delta B_{k/n} |$. Thus  we get
$$  \phi( \XE_{k/n},  \Delta B_{k/n} ) =  \XE_{k/n} +  \sigma ( \XE_{k/n}) \Delta B_{k/n}      + \frac{1}{2}    (\sigma   \sigma') (\phi( \XE_{k/n}, \zeta_{k}))    (\Delta B_{k/n} )^{2} $$  by the definition of $\phi$.
So we finally obtain 
\begin{eqnarray*}
 && \phi( \XE_{(k+1)/n},-B_{(k+1)/n})-  \phi( \XE_{k/n},-B_{k/n}) \\ && \qquad  =  
\left[ b( \XE_{k/n})  \Delta - \frac{1}{2}  (\sigma   \sigma') (\phi( \XE_{k/n}, \zeta_{k}))    (\Delta B_{k/n})^{2} \ \right]   \frac{\partial \phi}{\partial x_{1}} (\xi_{k}, -B_{(k+1)/n}).
\end{eqnarray*}
Since
$$ \frac{\partial \phi}{\partial x_{1}} (x_{1},x_{2})= \exp \left(\int_{0}^{x_{2}} \sigma'(\phi(x_{1},s)) \,ds  \right) $$
 by Lemma \ref{prop_h} (e), we have
$$  \left| \frac{\partial \phi}{\partial x_{1}} (\xi_{k}, -B_{(k+1)/n}) \right| \leq \exp ( c  \,  \|B \|_{\infty}). $$
Due to the assumptions,  the drift and diffusion coefficients satisfy a linear growth condition, i.e., 
$$ |b(x)| \leq c   (1+x), \quad  \quad \quad |\sigma(x)| \leq c   (1+x)$$
for $x \in \R$.
Hence we get
\begin{eqnarray*}
  &&  |\phi( \XE_{(k+1)/n},-B_{(k+1/n)})| \, \leq \,  |\phi( \XE_{k/n},-B_{k/n})|  + c \, \exp ( c \,   \|B \|_{\infty})     (1+ | \XE_{k/n}|)   \Delta  \\ && \qquad \qquad \qquad \qquad \qquad \qquad +   c  \,  \exp ( c \,   \|B \|_{\infty})   (1+ |\phi( \XE_{k/n}, \zeta_{k})|)    (\Delta B_{k/n})^{2}.
\end{eqnarray*}
Since by Lemma  \ref{prop_h} (a)
$$  \XE_{k/n}=\phi(\phi( \XE_{k/n},-B_{k/n}),B_{k/n}),   $$ 
and  $\phi(0,0)=0$,
we have by Lemma \ref{prop_h} (b) and (e)
$$
| \XE_{k/n}|\leq c \,  \exp ( c \,   \|B \|_{\infty} )   \left(   \ |\phi( \XE_{k/n},-B_{k/n})| +  \| B \|_{\infty} \right)
$$
and furthermore, since $|\zeta_{k}| \leq |\Delta B_{k/n}|$ 
\begin{eqnarray*}
&& |\phi( \XE_{k/n}, \zeta_{k})|  \leq  
c  \, \exp ( c \,  \|B \|_{\infty} )   \left( 
 | \XE_{k/n}| +  \|B \|_{\infty}  \right) \\ && \qquad \qquad \quad \, \, \, \,  
\leq c \,  \exp ( c \,   \|B \|_{\infty} )    \left(   \ |\phi( \XE_{k/n},-B_{k/n})| +  \| B \|_{\infty} \right).
\end{eqnarray*}
Together with
$$ |\Delta B_{k/n}| \leq \|B \|_{1/2}   \Delta^{1/2}, $$
this yields 
\begin{eqnarray*}
  &&  |\phi( \XE_{(k+1)/n},-B_{(k+1/n)})| \, \leq \,  |\phi( \XE_{k/n},-B_{k/n})|   \left[ 1  + c\,  \exp ( c \,   \|B \|_{\infty})   (1+ \|B \|_{1/2}^2)   \Delta \right]  \\ && \qquad \qquad \qquad \qquad \qquad \qquad \qquad \qquad \qquad \qquad \quad \, +   c \, \exp ( c \,   \|B \|_{\infty})   (1+ \|B \|_{1/2}^{2})   \Delta.  
\end{eqnarray*}
Setting $M= \|B \|_{\infty} + \|B \|_{1/2}^{2}$ it follows by Lemma \ref{rek}
\begin{eqnarray*} 
&&  |\phi( \XE_{k/n},-B_{k/n})| \leq  \sum_{k=1}^{n}   c  \, \exp ( c  \,  M )  \Delta    \prod_{j=1}^{n} (1+  c\,  \exp ( c \,  M ) \Delta) \\ && \qquad \qquad \qquad  \quad \, \, \, \leq  \exp( c \, \exp ( c \,  M)) .
\end{eqnarray*}
Thus with
$$  \XE_{k/n} = \phi(\phi( \XE_{k/n},-B_{k/n}),B_{k/n}), $$
 we get the estimate
\begin{eqnarray*} && |  \XE_{k/n}|  \leq \exp(c  \, \| B \|_{\infty})   |\phi( \XE_{k/n},-B_{k/n})| +  c \exp(c  \, \| B\|_{\infty})   \| B \|_{\infty}  \leq \exp( c \, \exp ( c \,M )).
\end{eqnarray*}
\begin{flushright} $\square$ \end{flushright}
\bigskip
Now we will state some Lemmas, which will 
 be needed to determine the asymptotic constant of the error of the Euler scheme. The following Lemma \ref{prop_f}
can be shown by straightforward calculations.

\medskip

\begin{lemma}\label{prop_f}
Denote
$$ f(x,y)= \exp \left( - \int_{0}^{x} \sigma'(\phi(y,s)) \, ds  \right )   b(\phi(y,x)), \qquad x, y \in \R.$$
Then we have $f \in \mathcal{C}^{1,2}$ and in particular
\begin{eqnarray*}
 f_{y}(x,y)   =
 b'(\phi(y,x)) -f(x,y)   \int_{0}^{x} \sigma''(\phi(y,s)) \frac{ \partial \phi}{\partial x_{1}}(y,s) \, ds, \qquad x, y \in \R.
\end{eqnarray*}
\end{lemma}
\bigskip

\begin{lemma}{\label{lem_const}}
We have $a.s.$
\begin{eqnarray*}
 && \exp \left( \int_{s}^{t} b'(X_{u}) \, d u + \int_{s}^{t} \sigma'(X_{u}) d^{-} B_{u} \right) = \\ && \qquad \qquad \qquad 
  \frac {\partial \phi}{ \partial x_{1}}(A_t,B_t) \frac{\partial \phi}{ \partial x_{1} }(X_s,-B_s)    \exp \left( \int_{s}^{t} f_{y}(B_u,A_u) \, du\right)  , \qquad 0 \leq s \leq t \leq 1. 
\end{eqnarray*}
\end{lemma}
{\bf Proof.}
By Lemma \ref{prop_h} (d) and (e) we have
\begin{eqnarray*}
&& \quad \frac{\partial \phi}{ \partial x_{1} }(A_t,B_t)  =  \exp \left( \int_{0}^{B_t} \sigma'(\phi(A_t,u ) )\, du \right), \\ &&  \frac{\partial \phi}{ \partial x_{1} }(X_s,-B_s)   =  \exp \left( -\int_{0}^{B_s} \sigma'(\phi(A_s,u ) )\, du \right).
\end{eqnarray*}
Using the notation
$$ g(x,y)= \int_{0}^{x} \sigma'(\phi(y,u)) \, du $$
we get by Lemma \ref{prop_f}
\begin{eqnarray*}
 &&  \frac {\partial \phi}{ \partial x_{1}}(A_t,B_t) \frac{\partial \phi}{ \partial x_{1} }(X_s,-B_s )   \exp \left( \int_{s}^{t} f_{y}(B_u,A_u) \, du\right)  \\ &&\qquad \qquad  = \exp \left( g(B_t,A_t)-g(B_s,A_s)+ \int_{s}^{t} f_{y}(B_u,A_u) \, du\right)
 \\ && \qquad \qquad  =\exp \left( \int_{s}^{t} b'(X_u) \, du  \right)    \exp \left( g(B_t,A_t)-g(B_s,A_s) \right) \\ && \qquad \qquad \qquad \qquad   \exp \left( - \int_{s}^{t} 
 \int_{0}^{B_u} \sigma''(\phi(A_u,\tau)) \frac{ \partial \phi}{ \partial x_{1} }
(A_u,\tau) \, d \tau  \, f(B_{u},A_{u}) d u\right).
\end{eqnarray*} 
Since
$$g_{x}(x,y)= \sigma'(\phi(y,x))  $$
and
$$g_{y}(x,y)=  \int_{0}^{x} \sigma''(\phi(y,s)) \frac{ \partial \phi}{\partial x_{1}} (y,s) \, ds $$
we have by the change of variable formula for Riemann-Stieltjes integrals, see e.g., 
\cite{ZK},
\begin{eqnarray*}
 && g(B_t,A_t)-g(B_s,A_s) =  \int_{s}^{t} \sigma'(\phi(A_u,B_u)) \, d^{-}B_u \\ && \qquad \qquad \qquad \qquad \qquad \qquad  \qquad  + \int_{s}^{t}  \int_{0}^{B_u} \sigma''(\phi(A_u,v)) \frac{ \partial \phi}{ \partial x_{1} }(A_u,v) \, dv \, d A_u,
\end{eqnarray*}
Since $A'_{t}=f(B_{t},A_{t})$ we finally get 
\begin{eqnarray*}
 &&  g(B_t,A_t)-g(B_s,A_s)- \int_{s}^{t}  \int_{0}^{B_u} \sigma''(\phi(A_u,v)) \frac{ \partial \phi}{\partial x_{1}} (A_u,v) \, dv \, f(B_{u},A_{u}) d u \\ && \qquad  = \int_{s}^{t} \sigma'(X_u) \, d^{-}B_u,
\end{eqnarray*} which shows the assertion.
\begin{flushright} $\square$ \end{flushright}
\bigskip
The next Lemma can be shown by a density argument.
\medskip
\begin{lemma}\label{dens_arg}
Let $g, h \in \mathcal{C}([0,1])$ and denote $\Delta h_{k/n}=  h((k+1)/n) - h(k/n)$ for $k= 0 \ldots, n-1$, $n \in \N$. If  
$$ \sup_{t \in [0,1]} \left| n^{2H-1}   \sum_{k=0}^{n-1}{\bf 1}_{[0,t]}(k/n)(\Delta h_{k/n})^2 -t \right| \, \longrightarrow \, 0$$ as $n \rightarrow \infty$,
then it follows
$$ \sup_{t \in [0,1]} \left| n^{2H-1}   \sum_{k=0}^{n-1}g(k/n){\bf 1}_{[0,t]}(k/n)(\Delta h_{k/n})^2 - \int_{0}^{t}g(s) \, ds  \right| \, \longrightarrow \, 0$$ as $n \rightarrow \infty$.
\end{lemma}
\bigskip
Now we finally prove Theorem 
\ref{thm1} and \ref{thm3}.
In the following we will denote by
$C$  random constants, which depend only on $b$, $\sigma$, their derivatives, $x_{0}$, $H$, $\|B\|_{\infty}$ and $\|B\|_{\lambda}$ with $\lambda < H$,  regardless of their  value.
We start with the proof of Theorem  \ref{thm3}.

\bigskip
\noindent
{\bf Proof of Theorem \ref{thm3}.} (1) We first establish a rough estimate for the pathwise error of the Euler scheme. For this, we follow the lines of the proof of Theorem 2.6. in \cite{talay}.
Set $$ \DE_{k}=\phi( \XE_{k/n},-B_{k/n}), \qquad k=0,\ldots, n $$ for $n \in \N$.
By a  Taylor expansion, the properties of $\phi$ and  Lemma \ref{bound_euler} we have
\begin{eqnarray*}
&&  \DE_{k+1} - \DE_{k}  =   \frac{ \partial \phi }{\partial x_{1}} ( \XE_{k/n},-B_{k/n})   ( \XE_{(k+1)/n}- \XE_{k/n} )  - \frac{ \partial \phi }{\partial x_{2}} ( \XE_{k/n},-B_{k/n}) \Delta B_{k/n}\\ && \qquad \qquad \qquad \qquad \quad +  \frac{1}{2} \frac{ \partial^{2} \phi }{\partial x_{1}^{2}} ( \XE_{k/n},-B_{k/n})    \sigma( \XE_{k/n})^{2}   (\Delta B_{k/n})^{2} 
\\ && \qquad \qquad \qquad \qquad  \quad +  \frac{1}{2} \frac{ \partial^{2} \phi }{\partial x_{2}^{2}} ( \XE_{k/n},-B_{k/n})    (\Delta B_{k/n})^{2} 
 \\ && \qquad \qquad \qquad \qquad  \quad -  \ \frac{ \partial^{2} \phi }{\partial x_{2} \partial x_{1}} ( \XE_{k/n},-B_{k/n})   \sigma( \XE_{k/n})    (\Delta B_{k/n})^{2}
+
R_{k}^{(1)}
\end{eqnarray*} with
\begin{eqnarray} {\label{develop_D}}
|R_{k}^{(1)}| \leq C   (   (\Delta B_{k/n})^{3} + \Delta\cdot\Delta B_{k/n}  + \Delta^2  ). \end{eqnarray}
Since
\begin{eqnarray*}
 && -  \frac{1}{2} \frac{ \partial \phi}{\partial x_{1}} (\XE_{k/n},-B_{k/n}) (\sigma \sigma')( \XE_{k/n})  \\ && \qquad =\frac{1}{2} \frac{ \partial^{2} \phi }{\partial x_{1}^{2}} ( \XE_{k/n},-B_{k/n})    \sigma( \XE_{k/n})^{2}    +   \frac{1}{2} \frac{ \partial^{2} \phi }{\partial x_{2}^{2}} ( \XE_{k/n},-B_{k/n})   -  \frac{ \partial^{2} \phi }{\partial x_{2} \partial x_{1}} ( \XE_{k/n},-B_{k/n})   \sigma( \XE_{k/n}) 
\end{eqnarray*}
by Lemma \ref{prop_h} (c), we have

$$ \DE_{k+1} =\DE_{k}+ b( \XE_{k/n})    \frac{ \partial \phi}{\partial x_{1}} (\XE_k,-B_{k/n})   \D + \widehat{Q}_{k} + R_{k}^{(1)}, $$
for $k =0, \ldots, n-1$, $n \in \N$,
with
$$\widehat{Q}_{k}=    - \frac{1}{2} (\sigma \sigma')( \XE_{k/n})   (\Delta B_{k/n} )^{2}   \frac{ \partial \phi}{\partial x_{1}} (\XE_k,-B_{k/n}), \qquad k=0, \ldots, n.$$ 
Since $ \XE_{k/n}=\phi(\DE_{k},B_{k/n})$ and using Lemma \ref{prop_h} (d) and (e)
we get 
\begin{eqnarray*}
\DE_{k+1} =\DE_{k}+ f(B_{k/n}, \DE_{k})   \Delta  + \widehat{Q}_{k}
 + R_{k}^{(1)} \end{eqnarray*}
for $k =0, \ldots, n-1$, $n \in \N$,
with the function $f$ given in Lemma \ref{prop_f}. 
Note that
$$  \sup_{n \in \N} \sup_{k=0, \ldots, n} |\DE_{k}| \leq \exp(c \,  \exp (c \,(\|B\|_{\infty}+  \|B \|_{1/2}^{2} ))),$$
as a consequence of Lemma \ref{bound_euler}. Now set  $$e_{k}= A_{k/n}- \DE_{k}, \qquad k=0, \ldots n.$$ 
We have $e_{0}=A_0 - \phi(x_{0},0)=0$ and 
$$
|e_{k+1}|  \leq |e_{k}|   (1+ C   \Delta) +|Q_{k}|+ |R_{k}^{(1)}| + 
\left |   \int_{k/n}^{(k+1)/n} f(B_\tau,A_\tau)- f(B_{k/n},A_{k/n}) \, d \tau  \right|.$$
Since 
\begin{eqnarray*} 
\left | \int_{k/n}^{(k+1)/n} f(B_\tau,A_\tau)- f(B_{k/n},A_{k/n}) \, d \tau  \right| \leq C  \int_{k/n}^{(k+1)/n} |B_\tau-
B_{k/n}| \, d \tau + 
C    \Delta^{2} \end{eqnarray*}
we can rewrite the above recursion as
$$|e_{k+1}|  \leq |e_{k}|   (1+ C   \Delta) + |Q_{k}| + |R_{k}^{(2)}| $$
with 
\begin{eqnarray} |R_{k}^{(2)}|  \leq C   \left( \|B \|_{H-\eps}^{3}   \Delta^{3H-3\eps} + \|B \|_{H-\eps}   \Delta^{H+1 -\eps} + \Delta^{2} \right) \leq C    \Delta^{H+1 -\eps}.  \label{est_R2} \end{eqnarray} 
Since also
$$  |\widehat{Q}_{k}|  \leq C     (\D B_{k/n})^{2} \leq C   \Delta^{2H-2\eps}, $$ we get by Lemma \ref{rek}
\begin{eqnarray}
 \label{pre_est_D} \max_{k=0, \ldots, n }|e_{k}| \leq \prod_{i=0}^{n-1}(1 + C   \Delta)    \sum_{j=0}^{n-1} (|\widehat{Q}_{j}|+ |R_{j}^{(2)}| )\leq C   \Delta^{2H-1-2\eps}. \end{eqnarray}
Moreover, we have
\begin{eqnarray}\label{pre_est_X}
 \max_{k=0, \ldots, n }|X_{k/n} -  \XE_{k/n}| \leq C   \Delta^{2H-1-2\eps},
\end{eqnarray}
due to 
$ X_t = \phi(A_t,B_{t})$, $t \in [0,1]$, and $ \XE_{k/n}=\phi(\DE_{k},B_{k/n})$, $k=0, \ldots, n$. \\

\noindent
(2) Now we derive the exact asymptotics of the error of the Euler scheme. We can write the recursion for the error $e_{k}=A_{k/n}-\DE_{k}$  as
\begin{eqnarray*}
 && e_{k+1} = e_{k}+ f_{y}(B_{k/n},A_{k/n})   e_{k}   \Delta + \widehat{Q}_{k} 
 + R_{k}^{(2)}  +  \frac{1}{2}   f_{yy}(B_{k/n},\eta_{k})   e_{k}^{2}   \Delta
\end{eqnarray*} with $\eta_{k}$ between $A_{k/n}$ and $\DE_{k}$. Put
$$ Q_{k}= -\frac{1}{2}(\sigma \sigma')(X_{k/n})   (\D B_{k/n})^{2}   \frac{\partial \phi}{\partial x_{1}} (X_{k/n}, -B_{k/n}). $$ 
By (\ref{pre_est_D}) we have
$$| Q_{k}-  \widehat{Q}_{k}| \leq C   \Delta^{ 4H-1-4 \eps}. $$
Since moreover
$$ 
\left| \int_{k/n}^{(k+1)/n} f_{y}(B_{t},A_{t}) \, dt - f_{y}(B_{k/n},A_{k/n})   \Delta \right| \leq C   \Delta^{H+1-\eps},$$
we get by   (\ref{est_R2}) and (\ref{pre_est_X}) 
\begin{eqnarray*} 
 e_{k+1} = e_{k}+    e_{k}   \int_{k/n}^{(k+1)/n} f_{y}(B_{t},A_{t}) \, dt  + Q_{k} +R_{k}^{(3)}
\end{eqnarray*} with
$$ |R^{(3)}_{k}| \leq C   \Delta^{\min \{ 4H-1-4 \eps, H +1-\eps \} }. $$
Applying Lemma \ref{rek} yields
\begin{eqnarray*}
A_{k/n}-\DE_{k}=  \sum_{i=0}^{k-1} Q_{i}     \prod_{j=i+1}^{k-1} \left (1+ \int_{j/n}^{(j+1)/n} f_{y}(B_{t},A_{t}) \, dt \right)     +   R_{k}^{(4)},
\end{eqnarray*}
with
\begin{eqnarray*}
 && \sup_{k=0, \ldots, n} |R_{k}^{(4)}|  = \sup_{k=1, \ldots, n}   \left|      \sum_{i=0}^{k-1} R^{(3)}_{i}       \prod_{j=i+1}^{k-1} \left (1+ \int_{j/n}^{(j+1)/n} f_{y}(B_{t},A_{t}) \, dt \right)  \right| \\ &&  \qquad \qquad \qquad \leq C   \Delta^{\min \{ 4H-2-4 \eps, H -\eps \} }. \end{eqnarray*}
Thus it remains to consider the term
$$  \sum_{i=0}^{k-1} Q_{i}    \prod_{j=i+1}^{k-1} \left (1+ \int_{j/n}^{(j+1)/n} f_{y}(B_{t},A_{t}) \, dt \right). $$
Now set 
$$ a_{j}= \int_{j/n}^{(j+1)/n} f_{y}(B_{t},A_{t}) \, dt, \qquad j=0, \ldots, n-1.$$
Since
\begin{eqnarray*}
&& \left| \prod_{j=i+1}^{k-1} (1+ a_{j}) - \exp\left( \sum_{j=i}^{k-1} a_{j}\right) \right| \\ && \qquad \leq  \exp \left(- \sum_{j=0}^{i-1} a_{j} \right) \left| \exp \left( \sum_{j=0}^{i-1} a_{j} \right) \prod_{j=i+1}^{k-1}  (1+ a_{j}) - \exp\left( \sum_{j=0}^{k-1} a_{j}\right) \right|
\end{eqnarray*}
we get 
$$ \sup_{0  \leq j < k-1 \leq n}\left| \prod_{j=i+1}^{k-1} (1+ a_{j}) - \exp\left( \sum_{j=i}^{k-1} a_{j}\right) \right|  \leq C \Delta$$
by a straightforward application of  Lemma \ref{rek}.
Hence we  have 
\begin{eqnarray}{\label{err_develop_D}}
A_{k/n}-\DE_{k}=  \sum_{i=0}^{k-1} Q_{i}    \exp \left( \int_{i/n}^{k/n} f_{y}(B_s,A_s) \, ds \right)  + R_{k}^{(5)}.
\end{eqnarray}
with
$$ \sup_{k=0, \ldots, n} |R_{k}^{(5)}|  \leq C   \Delta^{\min \{ 4H-2-4 \eps, H -\eps \} }. $$
Moreover, since
$X_t=\phi(A_{t},B_{t})$, $t \in [0,1]$  and $ \XE_{k/n}=\phi(\DE_{k},B_{k/n})$, $k=0, \ldots, n$ we have 
$$ X_{k/n}- \XE_{k/n}= \frac{\partial \phi}{ \partial x_{1}} (A_{k/n},B_{k/n})   (A_{k/n}-\DE_{k}) 
  + \frac{1}{2} \frac{\partial^{2} \phi}{ \partial x_{1}^{2}}(\theta_{k},B_{k/n}) (A_{k/n}-\DE_{k})^{2}$$
with $\theta_{k}$ between $A_{k/n}$ and $\DE_{k}$.
It follows by  (\ref{pre_est_D}) and (\ref{err_develop_D})
$$ X_{k/n}- \XE_{k/n}=  \frac{\partial \phi}{ \partial x_{1}} (A_{k/n},B_{k/n})   \sum_{i=0}^{k-1}   Q_{i}   
 \exp \left( \int_{i/n}^{k/n} f_{y}(B_s,A_s) \, ds \right)  + R_{k}^{(6)}$$
with
$$ |R_{k}^{(6)}| \leq C   \Delta^{\min \{ 4H-2-4 \eps, H -\eps \} }. $$
Since finally by Lemma \ref{lem_const}
\begin{eqnarray*}
&&  \frac{\partial \phi}{ \partial x_{1}} (A_{k/n},B_{k/n})   Q_{i} 
 \exp \left( \int_{i/n}^{k/n} f_{y}(B_s,A_s) \, ds \right)
 \\ && \qquad \qquad =   -\frac{1}{2} \frac{\partial \phi}{ \partial x_{1}} (A_{k/n},B_{k/n})  
 (\sigma \sigma')(X_{i/n})   \frac{\partial \phi}{\partial x_{1}} (X_{i/n}, -B_{i/n}) \\ && \qquad \qquad \qquad \qquad    \exp \left( \int_{i/n}^{k/n} f_{y}(B_s,A_s) \, ds \right)   (\D B_{i/n})^{2}  \\ && \qquad \qquad   = -\frac{1}{2} \sigma'(X_{i/n})   D_{i/n}X_{k/n}     (\D B_{i/n})^{2},
\end{eqnarray*}
we have 
\begin{eqnarray}{\label{err_develop_X}}
X_{k/n}- \XE_{k/n}= - \frac{1}{2} \sum_{i=0}^{k-1} \sigma'(X_{i/n})   D_{i/n}X_{k/n}     (\D B_{i/n})^{2}  + R_{k}^{(6)}.
\end{eqnarray}
Define 
$\widetilde{X}_{t}^{n}= X_{k/n}$
for $ t \in [k/n, (k+1)/n[$.
Since clearly
$$ \| X - \widetilde{X}^{n} \|_{\infty} \leq C    \Delta^{H-\eps}, $$
we have
$$  \| X -\XE \|_{\infty} = \max_{k=0, \ldots, n}|X_{k/n}- \XE_{k/n}| + R_{k}^{(7)}$$ with
$$ |R_{k}^{(7)}| \leq C   \Delta^{H-\eps}.$$
Thus it follows $$  \lim_{n \rightarrow \infty} \,\, n^{2H-1}   \| X -\XE \|_{\infty} \, =  \, 
\lim_{n \rightarrow \infty} \,\, \max_{k=0, \ldots, n} n^{2H-1}   |X_{k/n}- \XE_{k/n}| $$ 
and we get by (\ref{err_develop_X})
\begin{eqnarray*}    \lim_{n \rightarrow \infty} \,\, n^{2H-1}   \| X -\XE \|_{\infty} \, =   \,  \lim_{n \rightarrow \infty} \,\, \max_{k=0, \ldots, n}   \frac{n^{2H-1}}{2}   \left| \sum_{i=0}^{k-1} \sigma'(X_{i/n}) D_{i/n}X_{k/n}     (\D B_{i/n})^{2} \right|.  \end{eqnarray*}
Furthermore it holds
\begin{eqnarray*}
 &&  \lim_{n \rightarrow \infty} \,\, \max_{k=0, \ldots, n}   \frac{n^{2H-1}}{2}   \left| \sum_{i=0}^{k-1} \sigma'(X_{i/n}) D_{i/n}X_{k/n}     (\D B_{i/n})^{2} \right|
\\ && \qquad \qquad \qquad =  \lim_{n \rightarrow \infty} \,\, \sup_{t \in [0,1]} \,   \frac{n^{2H-1}}{2}   Z_{t}   \left| \sum_{i=0}^{n-1} {\bf 1}_{[0,t]}(i/n) \sigma'(X_{i/n}) D_{i/n}X_{1}     (\D B_{i/n})^{2} \right|.
\end{eqnarray*}
with
$$Z_{t}=\exp \left( -\int_{t}^{1} b'(X_{u}) \, du - \int_{t}^{1} \sigma'(X_{u}) \, d^{-} B_{u} \right), \qquad t \in [0,1].$$ This is due to the fact that the sample paths of $Z$ are H\"older continuous of any order $\lambda <H$.
It is well known that 
$$  n^{2H-1}   \sum_{k=0}^{n-1}{\bf 1}_{[0,t]}(k/n) (\Delta B_{k/n})^2 \stackrel{a.s.}{\longrightarrow} t$$ as $n \rightarrow \infty$ for all $t \in [0,1]$.
Since  $n^{2H-1}   \sum_{k=0}^{n-1}{\bf 1}_{[0,t]}(k/n) (\Delta B_{k/n})^2$ is monotone in $t$, the exceptional set of the almost sure convergence can be chosen independent of $t \in [0,1]$. Thus we get  by Dini's second theorem that a.s.
\begin{eqnarray} \label{dini}  \lim_{n \rightarrow \infty}  \sup_{t \in [0,1]} \left| n^{2H-1}   \sum_{k=0}^{n-1}{\bf 1}_{[0,t]}(k/n) (\Delta B_{k/n})^2  -t \right|=0. \end{eqnarray}
Hence it follows by Lemma \ref{dens_arg}
\begin{eqnarray*}
 &&    \lim_{n \rightarrow \infty} \,\, \sup_{t \in [0,1]} \,   n^{2H-1}  Z_{t}   \left| \sum_{i=0}^{n-1} {\bf 1}_{[0,t]}(i/n) \sigma'(X_{i/n}) D_{i/n}X_{1}  (\D B_{i/n})^{2} \right| 
\\ && \qquad \qquad  \qquad \quad = \,\,  \sup_{t \in [0,1]} \,   Z_{t} \left|\int_{0}^{t} \sigma'(X_{u}) D_{u}X_{1} \, du \right|  \qquad a.s.,
\end{eqnarray*} which finally shows the assertion.
\begin{flushright} $\square$ \end{flushright}

\bigskip
\bigskip
\noindent
{\bf Proof of Theorem \ref{thm1}.}
By (\ref{err_develop_X}) we have
$$ X_1-\XE_{n}= - \frac{1}{2} \sum_{i=0}^{n-1} \sigma'(X_{i/n})   D_{i/n}X_{1}     (\D B_{i/n})^{2}  + R_{n}^{(6)}$$
with
$$|R_{n}^{(6)}| \leq C    \Delta^{\min \{ 4H-2-4 \eps, H -\eps \} }.$$
The assertion follows then by (\ref{dini}) and Lemma \ref{dens_arg}, as in the previous proof.
\begin{flushright} $\square$ \end{flushright}

{\bf  Acknowledgements.} 
The authors are grateful to an anonymous referee for a careful and thorough reading of this work
and also for his valuable suggestions.

\end{document}